\documentstyle[12pt,a4]{article}

\newtheorem{lemma}{Lemma}[section]
\newtheorem{theorem}[lemma]{Theorem}
\newtheorem{prop}[lemma]{Proposition}
\newtheorem{cor}[lemma]{Corollary}
\begin{document}

\section{Introduction}

Much of the recent work on finite completely primary rings has demonstrated
the fundamental importance of these rings in the structure theory of finite
rings with identity. Let $R$ be a finite ring. It turns out that $R$ has a
unique maximal ideal if and only if it is a full matrix ring over a
completely primary ring. In particular, rings with a unique maximal ideal
are not necessarily completely primary. Therefore, the study of rings with a
unique maximal ideal (i.e. Local rings) reduces to the study of completely
primary rings.
\vspace{0.2in}

More evidence for the importance of completely primary rings comes from the
fact that any commutative ring is a direct sum of completely primary rings.
Moreover, any finite ring $R$ is of the form $S+N$, where $S\cap N=\{0\}$
with $N$ a subgroup of the Jacobson Radical of $R$ and $S$ a direct sum as
an additive abelian group of certain matrix rings over completely primary
rings (see [7]).
\vspace{0.2in}

In this paper we consider rings of characteristic $p$ with property(T)
 (see [1]). Clearly, such rings are completely primary. The rings of
 characteristics $p^2$ and $p^3$ with property(T) will be
 considered in later work.
\vspace{0.2in}

In Section 2, we collect some preliminary results on finite completely
primary rings. In Section 3, we give a construction of rings with
property(T) and characteristic $p,$ and in Section 4, we formulate the
isomorphism problem of these rings. Section 5 considers the problem of
enumerating certain cases of these rings.

\section{Preliminaries}

 Let $R$ be a finite ring. Then the following results will be assumed,
 and for details the reader is referred to [1], [4] and [9]:
 \vspace{0.2in}

 \noindent {\bf 2.1}\quad  Every  element in $R$ is either a zero-divisor or
a unit, and there is no distinction between left and right zero-divisors
(units).
\vspace{0.2in}

\noindent {\bf 2.2}\quad If $R$ is also  completely primary with
characteristic $p^k$ and Jacobson radical $M$, then 

\noindent \quad(i)\quad $\left| R\right| =p^{nr}$, for some positive integers
$n$ and $r$ such that $k\leq n$;

\noindent \quad(ii)\quad $R/M\cong GF(p^r),$ the field of $p^r$ elements; 

\noindent \quad(iii)\quad If $k=n$, then $R={\bf Z}_{p^k}[b]$, where $b$ is an
element of $R$ of multiplicative order $p^r-1;$
$M=pR$ and $Aut(R)\cong Aut(R/pR).$

\noindent The rings in 2.2(iii) shall be denoted by $GR(p^{nr},p^n)$ and are called
$Galois~ rings$.
\vspace{0.2in}

\noindent {\bf 2.3}\quad If $R$ is also completely primary with
characteristic $p^k$ and Jacobson radical $M$ such that
 $\left| R/M\right| =p^r$, then $R$ 
 has a coefficient subring $R_o$ of the form $GR(p^{kr},p^k)$ which is
clearly a maximal Galois subring of $R$. Furthermore, any two coefficient
subrings are conjugate in $R$ and 
there exist $m_1,\ldots ,m_h\in M$ and
$\sigma _1,\ldots ,\sigma _h\in Aut(R_o)$ such that

\noindent \quad(i)\quad $R=R_o\bigoplus \sum_{i=1}^h\bigoplus R_im_i$
(as $R_o$-modules);

\noindent \quad(ii)\quad $m_ir=r^{\sigma _i}m_i$, for every $r\in R_o$.

\noindent The  $\sigma _i$ are uniquely determined by $R$ and $R_o$ and are
 called  the automorphisms associated
with the $m_i$  with respect to $R_o.$

\section{Rings with property(T) and characteristic $p$}

 Let $F$ be the Galois field $GF(p^r)$.
For integers $s,t,\lambda $ with $%
1\leq t\leq s^2,$ $\lambda \geq 0,$ let $U,$ $V,$ $W$ be $s,$ $\lambda ,$ $%
t- $dimensional vector spaces over $F$, respectively. Since $F$ is
commutative, we can think of them as both left and right $F$-spaces. Let $%
\left( a_{ij}^k\right) $be $t$ compatible matrices of size $s\times s$ with
entries in $F,$ $\{\sigma _1,\ldots ,\sigma _s\},$ $\{\tau _1,\ldots ,\tau
_\lambda \},$ $\{\theta _1,\ldots ,\theta _t\}$ be sets of automorphisms of $%
F$ (with possible repetitions) and let $\{\sigma _i\}$ and $\{\theta _k\}$
satisfy the additional condition that if $a_{ij}^k\neq 0,$ for any $k$ with $%
1\leq k\leq t,$ then $\theta _k=\sigma _i\sigma _j.$
Let $R$ be the additive group direct sum

$$R=F\oplus U\oplus V\oplus W.$$

Then $R$ may be given a ring structure via an appropriate multiplication
(see e.g. [1]). The ring $R$ is said to be given by Construction A,
and the following results are proved in [1]: 

\begin{theorem}
The ring $R$ given by Construction A is a ring with property(T) and of
characteristic   $p$. Conversely, every ring with property(T) and
characteristic $p$ is isomorphic to one given by Construction A.
 \end{theorem}

\begin{theorem}
Let $R$ be a ring of Construction A. Then the field $F$ lies in the
centre of $R$ if and only if $\sigma_i = \tau_{\mu} = \theta_k =id_F$,
for all $i=1,\ldots ,s;~\mu =1,\ldots ,\lambda;~k=1,\ldots ,t$; and $R$ is
 commutative if and only if $a^{k}_{ij}=a^{k}_{ji}$, for all
 $i,j=1,\ldots ,s$.
 \end{theorem}

In what follows, the integers $p$, $n$, $r$, $s$, $t$, and $\lambda$,
shall be called the invariants of $R$.
\vspace{0.2in}

  It is clear that what we have named invariants are indeed that, that is,
  isomorphic rings have that same invariants. On the other hand, it is
  easy to find examples of non-isomorphic rings with property(T) and
  characteristic $p$ with the same invariants.

  \section{The isomorphism problem}

In this section, we formulate the isomorphism problem of rings with
property(T) and characteristic $p$. We know that all rings of this type
are rings of Construction A. So, since $M^2\subseteq ann(M)$, we can
write
$$R = F\oplus U\oplus N,~{\rm where}~ N = V\oplus W,$$
and if we define $v_1,\ldots ,v_\lambda$ by $w_{t+1},\ldots , w_{t+\lambda}$,
and $\tau_1,\ldots ,\tau_{\lambda}$ by $\theta_{t+1},\ldots ,
\theta_{t+\lambda}$, respectively, then the multiplication in $R$
becomes
$$(\alpha _o,~\sum_i\alpha _iu_i,~\sum_{k=1}^{t+\lambda}\gamma
_kw_k)\cdot (\alpha _o^{^{\prime }},~\sum_i\alpha _i^{^{\prime }}u_i,
~\sum_{k=1}^{t+\lambda}\gamma _k^{^{\prime }}w_k)$$

$$ = (\alpha_o\alpha _o^{^{\prime }},~\sum_i[\alpha _o\alpha _i^{^{\prime }}+\alpha
_i(\alpha _o^{^{\prime }})^{\sigma _i}]u_i,
~\sum_k[\alpha _o\gamma _k^{^{\prime }}+\gamma _k(\alpha _o^{^{\prime
}})^{\theta _k}+
\sum_{i,j=1}^sa_{ij}^k\alpha _i(\alpha _j^{^{\prime}})^{\sigma _i}]w_k).$$
where $a_{ij}^{k}=0$, for all $k=t+\lambda$, $\lambda\geq 1$.
\vspace{0.2in}

Let $R$ be the ring given by the above multiplication with respect to
the compatible matrices $(a_{ij}^{k})$, with entries from $F$, and
automorphisms $\sigma_i,~\theta_k\in Aut(F)$ $(i=1,\ldots ,s;~
k=1,\ldots , t+\lambda)$; with $\theta_k = \sigma_i\sigma_j$ for any
$k$ with $1\leq k\leq t$, if $a_{ij}^{k}\neq 0$. Let $A = \{(a_{ij}^{k})
: k=1,\ldots , t\}$, and let us denote the ring $R$ with the above
multiplication by $R(A,\sigma_i, \theta_k)$.
\vspace{0.2in}

Thus, up to isomorphism, the ring $R(A,\sigma_i, \theta_k)$ is given by
the $t$ compatible matrices $(a_{ij}^{k})$ and the automorphisms
$\sigma_i,~\theta_k$, where $\sigma_i$ and $\theta_k$ occur with
multiplicity $n_i$ and $n_k$, respectively $(i=1,\ldots , s;~k=t+1,\ldots ,
t+\lambda)$.
\vspace{0.2in}

Let now $R^{'}$ be another ring of the same type with the same invariants
$p$, $n$, $r$, $s$, $t$, $\lambda$;
$$R^{'}=F\oplus U^{'}\oplus N^{'}, ~{\rm where}~ N^{'}=V^{'}\oplus W^{'},$$
with respect to compatible matrices $D=\{(d_{ij}^{k}): k=1,\ldots , t\}$
and associated automorphims $\sigma_{i}^{'}, ~\theta_{k}^{'}$. Let
$\sigma_{i}^{'}$ and $\theta_{k}^{'}$ occur with multiplicity
$n_{i}^{'}$ and $n_{k}^{'}$, respectively, and denote $R^{'}$ by
$R(D, \sigma_{i}^{'}, \theta_{k}^{'})$.
\vspace{0.2in}

We assume that the rings $R(A,\sigma_i, \theta_k)$ and
$R(D,\sigma_{i}^{'}, \theta_{k}^{'})$ are constructed from a common
maximal Galois subfield $F$.
\vspace{0.2in}

We introduce the symbol $M^\sigma$ to denote
$\sigma((a_{ij}))$ if $M=(a_{ij})$.

\begin{lemma}
With the above notations,
$$R(A,\sigma_i, \theta_k)\cong R(D,\sigma_{i}^{'}, \theta_{k}^{'})$$
if and only if there exist $B=(\beta_{\rho k})\in GL(t, F)$,
$C\in GL(s, F)$, $\sigma \in Aut(F)$ such that
$$D_\rho =\sum_{k=1}^{t}\beta_{k\rho}C^{T}A_{k}^{\sigma}C^{\sigma_\mu};$$
$\{\sigma_1,\ldots ,\sigma_s\}=\{\sigma_{1}^{'},\ldots ,\sigma_{s}^{'}\}$,
$\{\theta_{t+1},\ldots ,\theta_{t+\lambda}\}=\{\theta_{t+1}^{'},\ldots ,
\theta_{t+\lambda}^{'}\}$ and (after possible reindexing),
$n_{i}=n_{i}^{'}$, $n_k =n_{k}^{'}$  for $i=1,\ldots ,s;~k=t+1,\ldots ,
t+\lambda$.
\end{lemma}

\noindent{\bf Proof} Suppose there is an isomorphism
$$\phi : R(A,\sigma_i, \theta_k)\to
R(D,\sigma_{i}^{'}, \theta_{k}^{'}).$$
Then, $\phi(F)$ is a maximal Galois subfield of
$R(D,\sigma_{i}^{'}, \theta_{k}^{'})$ so there exits an invertible
element $w\in R(D,\sigma_{i}^{'}, \theta_{k}^{'})$ such that
$w\phi(F)w^{-1}=F$.

Now, consider the map
$$\psi:   R(A,\sigma_i, \theta_k) \to  R(D,\sigma_{i}^{'}, \theta_{k}^{'})$$
defined by $$ r  \mapsto  w\psi(r)w^{-1}.$$

Then, clearly, $\psi$ is an isomorphism from $R(A,\sigma_i, \theta_k)$ to
$R(D,\sigma_{i}^{'}, \theta_{k}^{'})$ which sends $F$ to itself.

Also
$$\psi(0,~\sum_{i}\alpha_{i}u_i,~0)= (0,~\sum_{\nu}\sum_{i}\psi(\alpha_i)
\alpha_{\nu i}u_{\nu}^{'},~y^{'})~~(y^{'}\in N^{'});$$
and
$$\psi(0,~0,~\sum_{k}\gamma_k w_k)=(0,~0,~\sum_{\rho}\sum_{k}
\psi(\gamma_{k})\beta_{\rho k}w_{\rho}^{'}).$$

Therefore,
$$\psi(0,~\sum_{i}\alpha_{i}u_i,~0)\cdot \psi(0,~\sum_{i}\alpha_{i}^
{'}u_i,~0)$$
$$=(0,~\sum_{\nu}\sum_{i}\psi(\alpha_i)\alpha_{\nu i}u_{\nu}^{'},~y^{'})
\cdot  (0,~\sum_{\nu}\sum_{i}\psi(\alpha_{i}^{'})\alpha_{\nu i}u_{\nu}^{'},
~y^{''})$$
$$= (0,~0,~\sum_{\rho}\sum_{\nu, \mu =1}^{s}\sum_{i,j=1}^{s}\psi(\alpha_i)
\psi(\alpha_{j}^{'})^{\sigma_\nu}\alpha_{\nu i}\alpha_{\mu j}^{\sigma_\nu}
d_{\nu \mu}^{\rho}w_{\rho}^{'}).$$

On the other hand,
$$\psi((0,~\sum_{i}\alpha_{i}u_i,~0)\cdot (0,~\sum_{i}\alpha_{i}^
{'}u_i,~0)) = \psi(0,~0,~\sum_{k}\sum_{i,j=1}^{s}\alpha_{i}(\alpha_{j}^{'})^
{\sigma_{i}}a_{ij}^{k}w_k)$$
$$=(0,~0,~\sum_{\rho}\sum_{k=1}^{t}\sum_{i,j=1}^{s}\psi(\alpha_{i}
(\alpha_{j}^{'})^{\sigma_{i}})\beta_{\rho k}\psi(a_{ij}^{k})w_{\rho}^{'}).$$

It follows that
$$\sum_{\nu, \mu =1}^{s}\sum_{i,j=1}^{s}\psi(\alpha_i)\psi(\alpha_{j}^{'})
^{\sigma_\nu}\alpha_{\nu i}\alpha_{\mu j}^{\sigma_\nu}d_{\nu \mu}^{\rho}
=\sum_{k=1}^{t}\sum_{i,j=1}^{s}\psi(\alpha_{i}
(\alpha_{j}^{'})^{\sigma_{i}})\beta_{\rho k}\psi(a_{ij}^{k}).\qquad 4.1$$

Now, $\psi|_F$ is an automorphism of $F$, and therefore, $\psi(a_{ij}^{k})
=\sigma(a_{ij}^{k})$, for some $\sigma \in Aut(F)$. Hence, $\sigma_\nu =
\sigma_i$, for all $i, \nu =1,\ldots , s.$
Hence, equation 4.1 now implies that
$$E^{T}D_{\rho}E^{\sigma_{\mu}}=\sum_{k=1}^{s}\beta_{k\rho}A_{k}^{\sigma},~
with~E=(\alpha_{\mu j});$$
that is
$$D_{\rho}=C^{T}[\sum_{k=1}^{t}\beta_{k\rho}A_{k}^{\sigma}]C^{\sigma_{\mu}}
=\sum_{k=1}^{t}\beta_{k\rho}C^{T}A_{k}^{\sigma}C^{\sigma_{\mu}},$$
where $C=E^{-1}$, as required.

That $\{\sigma_1,\ldots ,\sigma_s\}=\{\sigma_{1}^{'},\ldots ,
\sigma_{s}^{'}\}$,
$\{\theta_{t+1},\ldots ,\theta_{t+\lambda}\}=\{\theta_{t+1}^{'},\ldots ,
\theta_{t+\lambda}^{'}\}$ and (after possible reindexing),
$n_{i}=n_{i}^{'}$, $n_k =n_{k}^{'}$  for $i=1,\ldots ,s;~k=t+1,\ldots ,
t+\lambda$; follows from the fact that $R(A, \sigma_i, \theta_k)$ and
$R(D, \sigma_{i}^{'}, \theta_{k}^{'})$ are constructed from a common
maximal Galois subfield $F$.
\vspace{0.2in}

Now, suppose that there exist $B=(\beta_{\rho k})\in GL(t, F)$,
$C\in GL(s, F)$, $\sigma \in Aut(F)$ such that
$$D_\rho =\sum_{k=1}^{t}\beta_{k\rho}C^{T}A_{k}^{\sigma}C^{\sigma_\mu};$$
with $\{\sigma_1,\ldots ,\sigma_s\}=\{\sigma_{1}^{'},\ldots ,\sigma_{s}^{'}\}$,
$\{\theta_{t+1},\ldots ,\theta_{t+\lambda}\}=\{\theta_{t+1}^{'},\ldots ,
\theta_{t+\lambda}^{'}\}$ and (after possible reindexing),
$n_{i}=n_{i}^{'}$, $n_k =n_{k}^{'}$  for $i=1,\ldots ,s;~k=t+1,\ldots ,
t+\lambda$.

Consider the map
$$\psi:  R(A, \sigma_{i}, \theta_{k})  \to  R(D, \sigma_{i}^{'},
\theta_{k}^{'})$$
given by
 $$(\alpha_o,~\sum_{i}\alpha_{i}u_{i},~\sum_{k}\gamma_{k}w_{k}) \mapsto 
 (\alpha_{o}^{\sigma},~\sum_{\nu}\sum_{i}\alpha_{i}^{\sigma}\alpha_{\nu i}
 u_{\nu}^{'},~\sum_{\rho}\sum_{k}\gamma_{k}^{\sigma}\beta_{k\rho}
 w_{\rho}^{'}).$$

Then, it is easy to verify that $\psi$ is an isomorphism of the ring
$R(A, \sigma_{i}, \theta_{k})$ onto the ring
$R(D, \sigma_{i}^{'}, \theta_{k}^{'})$.

\begin{cor}
Let $A$ and $D$ be sets of compatible matrices with entries from $F$. If
$A$ and $D$ generate the same vector space over $F$,and if
$\sigma_i = \sigma_{i}^{'}$, $\theta_k = \theta_{k}^{'}$ with
$n_i =n_{i}^{'}$, $n_k =n_{k}^{'}$, then
$$R(A,\sigma_i, \theta_k)\cong R(D,\sigma_{i}^{'}, \theta_{k}^{'}).$$
\end{cor}

\section{The Enumeration problem}
In this section, we consider the problem of finding the number of
distinct (up to isomorphism) types of rings with property(T) and
characteristic $p$.
We find those rings of Construction A which give rise to
distinct non-isomorphic rings.
\vspace{0.2in}

We consider this for certain cases.

\subsection{The case where $s=1$}
For this case, $R$ is a ring of
Construction A with $t=1$, $\lambda\geq 0.$ Then, the only parameters in
the definition of $R$ are the automorphisms $\sigma_1$, $\theta_k$,
$(k=1+\lambda,~\lambda\geq 0)$, $\theta_1 =\sigma_{1}^{2}$; and the
element $a_{11}^{1}\in F^{*}$.
\vspace{0.2in}

Let us denote the ring $R$ by $R(a_{11}^{1}, \sigma_1, \theta_k)$. Thus,
up to isomorphism, the ring $R(a_{11}^{1}, \sigma_1, \theta_k)$ is
given by the element $a_{11}^{1}\in F^{*}$ and the
automorphisms $\sigma_1$, $\theta_k$, where for $\theta_k$,
$k>1$, $\theta_k$ occurs with multiplicity $n_k$. 
\vspace{0.2in}

If $R(d_{11}^{1}, \sigma_{1}^{'}, \theta_{k}^{'})$ is another ring of the
same type with the same invariants $p$, $n$, $r$, $s$, $t$, $\lambda$, with
$s=t=1$, then by Lemma 4.1

$$R(a_{11}^{1}, \sigma_1, \theta_k)\cong R(d_{11}^{1},
\sigma_{1}^{'}, \theta_{k}^{'})$$
if and only if there exist $\beta_{11}, ~\gamma \in F^{*}$ and
$\theta \in Aut(F)$ such that
$$d_{11}^{1}=\gamma \gamma^{\sigma_{1}}\beta_{11}(a_{11}^{1})^{\theta};~
\sigma_1 = \sigma_{1}^{'},~ \{\theta_2,\ldots , \theta_{1+\lambda}\} =
\{\theta_{1}^{'},\ldots , \theta_{1+\lambda}^{'}\}$$
and (after possible reindexing) $n_k = n_{k}^{'}$,  for every $k=2,\ldots ,
1+\lambda$.
\vspace{0.2in}

As a result of Lemma 4.1, if $\gamma,~\beta_{11}\in F^{*}$, then the rings
$R(a_{11}^{1}, \sigma_{1}, \theta_{k})$ and
$R(\gamma \gamma^{\sigma_{1}}\beta_{11}(a_{11}^{1})^{\theta}, \sigma_{1},
\theta_{k})$ are isomorphic. Hence, we can select $\gamma =1$ and
$\beta_{11} = ((a_{11}^{1})^{\theta})^{-1}$ to see that the rings
$R(a_{11}^{1}, \sigma_{1}, \theta_{k})$ and
$R(1, \sigma_{1}, \theta_{k})$ are isomorphic. So, counting the
isomorphism classes of the rings
$R(a_{11}^{1}, \sigma_{1}, \theta_{k})$ is a question of
counting the number of distinct ways of selecting the automorphisms.
\vspace{0.2in}

Consider now the automorphisms $\sigma_{1},~\theta_{2},\ldots ,
\theta_{1+\lambda}.$ Since $|Aut(F)|=r$,
the number of ways in which we can select $\sigma_1$
from $Aut(F)$ is $r$. Also, the number of ways we can select $\theta_2,\ldots
 ,\theta_{1+\lambda}$ from $Aut(F)$ ($\theta_k$ not necessarily distinct),
 is the number of solutions in the equation
 $$x_1 + x_2 + \ldots + x_r = \lambda$$
 in non-negative integers $x_1,~x_2, \ldots , x_r \in \{0, 1, \ldots ,
 \lambda \}$. This is well known to be (see [6], page 2)
 \begin{displaymath}
 \left( \begin{array}{c}
 r+\lambda -1\\
 \lambda
 \end{array} \right).
 \end{displaymath}
\vspace{0.2in}

Therefore, for a fixed $a_{11}^{1}\in F^{*}$, a $\sigma_{1}\in Aut(F)$
and a $\lambda$-selection of $\theta \in Aut(F)$, there is only one
ring up to isomorphism. Therefore, the number of isomorphism classes of
rings of Construction A of the same characteristic $p$ and same order,
with the same invariants $p$, $n$, $r$, $s$, $t$, $\lambda$, where
$s=t=1$, is
 \begin{displaymath}
 r\cdot \left( \begin{array}{c}
 r+\lambda -1\\
 \lambda
 \end{array} \right).
 \end{displaymath}
\vspace{0.2in}

We have thus proved the following

\begin{lemma}
The number of mutually non-isomorphic rings with property(T) and
characteristic $p$ and of the same order with the same invariants $p$, $n$,
$r$, $s$, $t$, $\lambda$, in which $s=t=1$, is
 \begin{displaymath}
 r\cdot \left( \begin{array}{c}
 r+\lambda -1\\
 \lambda
 \end{array} \right).
 \end{displaymath}
 Of these, only one is commutative, the others are not.
\end{lemma}

 If, in particular, $\lambda =0$, then the rings are principal ideal rings,
  so that we have 

 \begin{cor}
 The number of mutually non-isomorphic principal ideal rings with
 property(T) and characteristic $p$ (and of the same order) with the
 same invariants $p$, $n$, $r$ is $r$. Further, only one is commutative,
 the others are not.
 \end{cor}

 \subsection{The case where $t=s^{2}$}
Let $R$ be a ring of Construction A with the invariants $p$, $n$, $r$,
$s$, $t$, $\lambda$, where $t=s^2$; and let $\sigma_{i}$, $\theta_k$
 ($i=1, \ldots , s;~k=1, \ldots , t, t+1, \ldots , t+\lambda$) be the
associated automorphisms of $R$ with respect to a fixed maximal Galois
subfield $F$ of $R$, and let $A_1, A_2, \ldots , A_t$ be the compatible
structural matrices of $R$. Let ${\cal A}$ denote the subspace of
$M_{s}(F)$ generated by the matrices $A_1, \ldots , A_t$ over $F$.
Since $t=s^2$, ${\cal A} = M_{s}(F)$.
\vspace{0.2in}

Let now $R^{'}$ be another ring of the same type with the same invariants
$p$, $n$, $r$, $s$, $t$, $\lambda$, where $t=s^2$, with respect to the
automorphisms $\sigma_{j}^{'}, \theta_{l}^{'} \in Aut(F)$ ($j=1, \ldots ,
s;~l=1, \ldots, t, t+1, \ldots , t+\lambda$) and compatible structural
matrices $D_1, \ldots , D_t$, with respect to a common fixed maximal
Galois field $F$. Let ${\cal D}$ denote the subspace of $M_{s}(F)$
generated by $D_1, \ldots , D_t$ over $F$.  As before, since $t=s^2$,
then ${\cal D}=M_{s}(F)$.
\vspace{0.2in}

But ${\cal A} = {\cal D}=M_{s}(F)$. Thus, up to isomorphism, the rings
$R$ and $R^{'}$ are determined by the automorphisms $\sigma_i$, $\theta_k$
and $\sigma_{j}^{'}$, $\theta_{l}^{'}$ ($i,j=1, \ldots, s;~ k,l=t+1, \ldots ,
t+\lambda$), respectively.

\begin{lemma}
The number of mutually non-isomorphic rings with property(T) and
characteristic $p$, with the same invariants $p$, $n$, $r$, $s$, $t$,
$\lambda$, where $t=s^2$ is

$$ \left(
\begin{array}{c}
r + s - 1\\
s
\end{array}\right)\cdot
\left(
\begin{array}{c}
r + \lambda -1\\
\lambda
\end{array}\right).
$$
All of these rings are non-commutative.
\end{lemma}

\noindent{\bf Proof}  From the above discussion, it is clear that the number
of isomorphism classes in question is the number of ways in which we can
select $\{ \sigma_1, \ldots , \sigma_s\}$ and $\{ \theta_{t+1}, \ldots ,
\theta_{t+\lambda}\}$ ($\sigma_{i}, \theta_{k}$ not necessarily distinct)
from $Aut(F)$. Since $|Aut(F)|=r,$ the number in question is the
number of solutions of the two equations
$$ x_1 + x_2 + \ldots + x_r = s,$$
and
$$ y_1 + y_2 + \ldots + y_r = \lambda$$
in non-negative integers $x_1, \ldots , x_r \in \{0, 1,\ldots , s\}$
and $y_1, \ldots , y_r \in \{0, 1, \ldots , \lambda \}.$ This is well known
to be
$$ \left(
\begin{array}{c}
r + s - 1\\
s
\end{array}\right)\cdot
\left(
\begin{array}{c}
r + \lambda -1\\
\lambda
\end{array}\right).
$$

\subsection{The case where $F$ lies in the centre of $R$}
We now consider the case where the maximal Galois subfield
$F$ lies in the centre of $R$, that is, the case where all the associated
automorphisms of $R$ are equal to the identity automorphism (Theorem 3.2).
\vspace{0.2in}

We note that the description of the rings of this type reduces to the case
where $ann(M)$ coincides with $M^2$. Therefore, to enumerate the rings of this
type of a given order, say $p^{nr}$, where $ann(M)$ does not coincide with
$M^2$, we shall first write all the rings of this type of order $\leq p^{nr}$,
where $ann(M)$ coincides with $M^2$.
\vspace{0.2in}

In what follows, we assume that $ann(M)=M^2$.

\subsubsection{The case with $t=1$}
Suppose now that $R$ is a ring with property(T) and characteristic $p$ with
the invariants $p$, $n$, $r$, $s$, $t$; where $t=1$. Then, the ring $R$ is
defined by one structural matrix $A_1$, where $A_1$ is a non-zero $s\times s$
compatible matrix with entries from $F$.
\bigskip

Now, let $R(D_1)$ be another ring with property(T) and characteristic $p$
with the same invariants $p$, $n$, $r$, $s$, $t$, where $t=1$, and of
the same order as $R(A_1)$ and assume that they are constructed from a
common maximal Galois subfield $F$. Then, by Lemma 4.1, $R(A_1)\cong
R(D_1)$ if and only if there exists a $\sigma \in Aut(F)$, an
invertible matrix $C\in M_{s}(F)$ and a non-zero element $\beta \in F$
such that
$$D_1 = \beta^{-1}C^{T}A_{1}^{\sigma}C.$$

Congruence of matrices in the classical sense implies equivalence in the
sense defined above but not vice-versa as the following example shows:
\vspace{0.2in}

\noindent{\bf Example}  Let $F = {\bf F}_{4} =\{0, 1, \alpha, 1+\alpha \}$
and $\sigma \in Aut(F)$ such that $\sigma :x\mapsto x^2$, for every
$x\in F$. Consider the matrices
$$
\left(
\begin{array}{cc}
1 & 0\\
\alpha & 1
\end{array}\right),\qquad
\left(
\begin{array}{cc}
1 & 0\\
1+\alpha & 1
\end{array}\right) \in M_{2}(F).
$$
The $F$-spaces generated by these two matrices are equivalent since,
for instance,
$$
\left(
\begin{array}{cc}
1 & 0\\
\alpha & 1
\end{array}\right)^{\sigma}=
\left(
\begin{array}{cc}
1 & 0\\
1+\alpha & 1
\end{array}\right);
$$
while the two matrices are not congruent.
\vspace{0.2in}

However, in the cases where the automorphisms $\sigma$ can be reduced to
the identity (for instance, if $R$ is commutative or if $F$ is a prime field)
then equivalence comes very close to congruence (the element $\beta$ makes
the only difference). So, it makes sense to look at congruence classes.
\vspace{0.2in}

Let $N(s)$ denote the number of congruence classes of $s\times s$ matrices
over $F\cong GF(p^r)$. In [10], Waterhouse implicitly computed the number
of congruence classes of $n\times n$ matrices over finite fields, and in
[8], Newman obtained the number and representatives of congruence classes
of $n\times n$ symmetric matrices of positive rank $\leq n$ over
finite fields, and we restate these results here in our notation for easy
reference.

\begin{theorem}
$N(s)$ is the coefficient of $t^s$ in
$$\prod_{k\geq 1}(1 + t^k)^e(1 - qt^{2k})^{-1}(1 - t^k)^{-1},$$
where $e=1$ for even $q$ and $e=2$ for odd $q$.
\end{theorem}

\begin{theorem}

(i) Let $F$ be a finite field of characteristic $2$. Then every
symmetric matrix of $M_{n}(F)$ of rank $r$ is congruent to $I_{r}\oplus 0$
($r$ odd), or to $I_{r}\oplus 0$ or $(r/2)T\oplus 0$ ($r$ even), where
$kT$ denotes the direct sum of $k$ copies of
$$T=
\left(
\begin{array}{cc}
0 & 1\\
1 & 0
\end{array}
\right);
$$
and these are not congruent.

(ii) Let $F$ be a finite field of characteristic different from $2$. Then,
every symmetric matrix of $M_{n}(F)$ of rank $r$ is congruent to
$I_{r}\oplus 0$ or to $gI_{1}\oplus I_{r-1}\oplus 0$, where $g$ is a
fixed non-square in $F$; these are not congruent. Thus, the symmetric
matrices of any given rank fall into precisely two congruence classes.
\end{theorem}

We now consider the problem of finding the number of isomorphism classes of
rings with property(T) and characteristic $p$ with same invariants
$p$, $n$, $r$, $s$, $t$, $\lambda$; where $s>1$ and $t=1$. The solution of
 this problem depends on the much more difficult classical
problem of the classification of bilinear forms over finite fields.
\vspace{0.2in}

Consider the matrices $\beta^{-1}C^{T}AC$, where $A\in M_{s}(F_q)$.
\vspace{0.2in}

\noindent{\bf Case 1.} Suppose that $s=2$ and $t=1$. In [2], Bremser
obtained the congruence classes of matrices in $GL_2(F_q)$, and showed
that there are $q+3$ for odd $q$ and $q+1$ for even $q$. The number of
congruence classes in $M_2(F)$ over a finite field $F$ of any characteristic
$p$ can be calculated from the formula in Theorem 5.4 (see also
Waterhouse [10]), and here we give a complete set of representatives of these
classes, which include those obtained by Bremser in [2].
\vspace{0.2in}

\noindent{\bf (i)} $Char F \neq 2.$
\[
\left( \begin{array}{cc}
0 & 0\\
0 & 0 
\end{array}
\right), \quad \left( \begin{array}{cc}
0  & 1 \\
-1 & 0 
\end{array}
\right), \quad \left( \begin{array}{cc}
1  & 0 \\
0  & 0 
\end{array}
\right), \quad \left( \begin{array}{cc} 
1 & 0 \\ 
1 & 0  
\end{array} 
\right),
\]
\[
\left( \begin{array}{cc} 
g  & 0 \\ 
0  & 0  
\end{array} 
\right), \quad \left( \begin{array}{cc}  
g  & 0 \\  
2g & g   
\end{array}  
\right),\quad  
 \left( \begin{array}{cc}  
1  & 0 \\  
0  & 1   
\end{array}
\right), \quad  
 \left( \begin{array}{cc} 
1  & 0 \\ 
0  & g  
\end{array} 
\right), 
\]
\[
 \left( \begin{array}{cc}  
1      & 0 \\  
\gamma & 1   
\end{array}  
\right),\quad  
\left( \begin{array}{cc}  
    1  & 0 \\  
\gamma & g   
\end{array}  
\right),  
\]
where $\gamma$ runs over a complete set of coset representatives of
$\{ \pm 1\}$  in $F^{*}$; these are $q+7$ altogether. 
\vspace{0.2in} 

\noindent{\bf (ii)} $CharF=2$.

\[
\left( \begin{array}{cc} 
0  & 0 \\  
0  & 0   
\end{array}  
\right), \quad \left( \begin{array}{cc}  
1  & 0 \\  
0  & 0   
\end{array}  
\right), \quad \left( \begin{array}{cc}   
1  & 0 \\   
0  & 1    
\end{array}   
\right), \quad \left( \begin{array}{cc}   
0  & 1 \\   
1  & 0    
\end{array} 
\right), 
\]
\[ 
 \left( \begin{array}{cc}  
1  & 0 \\  
1  & 0   
\end{array}  
\right), \quad  \left( \begin{array}{cc}   
 1      & 0 \\
\alpha  & 1   
\end{array}
\right), \quad \alpha \in F^{*}  
\]      
and these are $q+4$ in all. 
\vspace{0.2in}

\noindent Now, suppose $|F|=2$. Then the non-zero congruence classes are   
\[
\left( \begin{array}{cc}
1  & 0 \\
0  & 0
\end{array}
\right), \quad \left( \begin{array}{cc}
1  & 0 \\
0  & 1
\end{array}
\right), \quad \left( \begin{array}{cc}
0  & 1 \\
1  & 0
\end{array}
\right),
\]
\[
 \left( \begin{array}{cc}
1  & 0 \\
1  & 0
\end{array}
\right), \quad  \left( \begin{array}{cc}
 1      & 0 \\
 1      & 1
\end{array}
\right).  
\]      
Since $\beta =1$ in this case, these matrices also represent  
 equivalence classes. Notice also that the equivalence class 
\[
\left( \begin{array}{cc}
1  & 0 \\
0  & 0  
\end{array}
\right) 
~{\rm contains~ the~ compatible~ matrix~}~ 
\left( \begin{array}{cc}
1  & 1 \\ 
1  & 1   
\end{array} 
\right);    
\] 
and therefore, we include this class among the equivalence classes that  
correspond to rings with property(T). 
Hence, the number of mutually non-isomorphic rings 
of this type is $5$, which is the number of non-zero congruence classes. 
\vspace{0.2in}

\noindent Suppose $|F|=p$, $p\neq 2$, then it can be deduced from the  
class representatives above that 
the number of non-zero congruence classes is $p+6$. Now, if $\beta =g$  
is an element of ${\bf F}^{*}_{p}$, it is easy to see that the congruence class 
\[
\left( \begin{array}{cc}
g   & 0 \\
2g  & g 
\end{array}
\right)
~{\rm is~ equivalent~to ~one~of~ the~classes ~of~the~form~}~
\left( \begin{array}{cc}
1       & 0 \\
\gamma  & 1
\end{array}
\right)
\] 
in (i) above. Also, the classes 
\[
\left( \begin{array}{cc}
1  & 0 \\
0  & 0 
\end{array}
\right)
~{\rm and~}~
\left( \begin{array}{cc}
g  & 0 \\
0  & 0 
\end{array}
\right)
\]
are equivalent. Moreover,  all the  equivalence 
classes  contain at least one compatible matrix.   
Therefore, the number of  equivalence classes in this case is $p+4$ and this 
also gives  the number and models for the corresponding rings. 
\vspace{0.2in}

\noindent {\bf Case 2.}\quad Suppose $s=3$ and $t=1$. In [5], 
B. Corbas and G. D. Williams have obtained the matrix representatives
for bilinear
forms on a three dimensional vector space over a finite field of any 
characteristic, without assuming that the form is symmetric or non-degenerate. 
We  give here a 
full list of the congruence classes as given in their main Theorem. 
\vspace{0.2in} 

\noindent {\bf (i)} $CharF\neq 2$. 
\[
\left( \begin{array}{ccc} 
0 & 0 & 0 \\
0 & 0 & 0 \\
0 & 0 & 0 
\end{array} 
\right), ~
\left( \begin{array}{ccc} 
1 & 0 & 0 \\ 
0 & 0 & 0 \\
0 & 0 & 0 
\end{array} 
\right), ~ 
\left( \begin{array}{ccc}  
\varepsilon  & 0 & 0 \\  
0            & 0 & 0 \\ 
0            & 0 & 0  
\end{array}  
\right), ~ 
\left( \begin{array}{ccc}  
1 & 0 & 0 \\  
0 & 1 & 0 \\ 
0 & 0 & 0  
\end{array}  
\right), 
\]
\[
\left( \begin{array}{ccc}  
1 &           0 & 0 \\  
0 & \varepsilon & 0 \\ 
0 &           0 & 0  
\end{array}  
\right), ~  
\left( \begin{array}{ccc}   
1  & 0 & 0 \\   
0  & 1 & 0 \\  
0  & 0 & 1   
\end{array}   
\right), ~  
\left( \begin{array}{ccc}   
1 & 0 &           0 \\   
0 & 1 &           0 \\  
0 & 0 & \varepsilon    
\end{array}   
\right),~ 
\left( \begin{array}{ccc}    
\mu  & 0  & 0 \\    
0    & 0  & 1 \\   
0    & -1 & 0    
\end{array}    
\right), 
\]
\[
\left( \begin{array}{ccc}   
\mu  &  0 & 0 \\   
0    &  0 & 0 \\  
0    &  1 & 0    
\end{array}   
\right), ~   
\left( \begin{array}{ccc}    
\mu  & 0            &           0 \\    
0    & \varepsilon  &           0 \\   
0    & 2\varepsilon & \varepsilon     
\end{array}    
\right), ~   
\left( \begin{array}{ccc}    
\mu & 0     & 0 \\    
0   & 1     & 0 \\   
0   & \gamma& 1      
\end{array}      
\right),~      
\left( \begin{array}{ccc}    
\mu  & 0      & 0 \\     
0    & 1      & 0 \\    
0    & \gamma & \varepsilon      
\end{array}     
\right),  
\]
\[
\left( \begin{array}{ccc}     
\mu  & 0 & 0 \\     
0    & 0 & 1 \\     
1    & 1 & 0       
\end{array}     
\right), ~    
\]
where $\mu \in \{ 0,~1,~\varepsilon \}$, with $\varepsilon$ an arbitrary but 
fixed non-square in $F^{*}$, and $\gamma$ runs over a complete set of coset 
representatives of $\{ 1,~-1 \}$ in $F^{*}$. Their total number is 
$3q+16$. 
\vspace{0.2in} 

\noindent {\bf (ii)} $CharF=2$. 
\[
\left( \begin{array}{ccc}
0  & 0 & 0 \\
0  & 0 & 0 \\
0  & 0 & 0   
\end{array}   
\right), ~
\left( \begin{array}{ccc} 
1  & 0 & 0 \\ 
0  & 0 & 0 \\ 
0  & 0 & 0   
\end{array}    
\right), ~
\left( \begin{array}{ccc} 
1  & 0 & 0 \\ 
0  & 1 & 0 \\ 
0  & 0 & 0   
\end{array}    
\right), ~ 
\left( \begin{array}{ccc}  
1  & 0 & 0 \\  
0  & 1 & 0 \\  
0  & 0 & 1   
\end{array}     
\right), 
\]
\[
\left( \begin{array}{ccc} 
0  & 0 & 0 \\ 
0  & 0 & 1 \\ 
0  & 1 & 0   
\end{array}    
\right), ~ 
\left( \begin{array}{ccc}  
\mu  & 0 & 0 \\  
0    & 0 & 0 \\  
0    & 1 & 0   
\end{array}     
\right), ~
\left( \begin{array}{ccc}  
\mu  & 0     & 0 \\  
0    & 1     & 0 \\  
0    &\gamma & 1    
\end{array}     
\right), ~  
\left( \begin{array}{ccc}   
1  & 0 & 0 \\   
0  & 0 & 0 \\   
1  & 1 & 0    
\end{array}      
\right),  
\] 
\[
\left( \begin{array}{ccc}   
1  & 0 & 0 \\   
0  & 0 & 1 \\   
1  & 1 & 0    
\end{array}      
\right), ~ 
\left( \begin{array}{ccc}   
1      & 0 & 0 \\   
0      & 0 & 1 \\   
\alpha & 1 & 1     
\end{array}      
\right), ~   
\]  
where $\mu \in \{ 0,~1 \}$, $\gamma \in F^{*}$ and $X^2 + \alpha X + 1$ is 
an arbitrary but fixed irreducible polynomial of degree two over $F$. Their 
total number is $2q+8$. 
\vspace{0.2in} 

\noindent Now, suppose $|F|=2$, then it can be deduced from the list of
class representatives in  (ii) above that the number of
non-zero congruence classes is $11$. Since
in this case $\beta = 1$, these classes also represent the equivalence 
classes of $1$-dimensional spaces of bilinear forms over ${\bf F}_2$. Further, 
since all the equivalence classes contain at least one compatible matrix, 
we conclude that the number of non-isomorphic rings with property(T) and 
characteristic $2$ with same invariants and with maximal Galois subfield 
${\bf F}_2$ is 11. The models of each of these is given by the corresponding 
equivalence class. 
\vspace{0.2in} 

\noindent Now, suppose $|F|=p$, $p\neq 2$. Then the list of class
representatives in (i) above gives
$3p+15$ non-zero congruence classes. As $\beta$ runs over 
the elements of ${\bf F}^{*}_{p}$, 
the congruence classes 
\[
\left( \begin{array}{ccc}
\varepsilon  & 0 & 0 \\
0            & 0 & 0 \\
0            & 0 & 0
\end{array}
\right), ~
\left( \begin{array}{ccc}
1 & 0 &           0 \\
0 & 1 &           0 \\
0 & 0 & \varepsilon
\end{array}
\right),~
\left( \begin{array}{ccc}
0    & 0            &           0 \\
0    & \varepsilon  &           0 \\
0    & 2\varepsilon & \varepsilon
\end{array}
\right), ~
\left( \begin{array}{ccc}
1    & 0            &           0 \\
0    & \varepsilon  &           0 \\
0    & 2\varepsilon & \varepsilon
\end{array}
\right),
\]
and 
\[
 \left( \begin{array}{ccc}
\varepsilon  & 0            &           0 \\
0            & \varepsilon  &           0 \\
0            & 2\varepsilon & \varepsilon
\end{array}
\right),~{\rm become~equivalent~to~} 
~\left( \begin{array}{ccc}
         1   & 0 & 0 \\
0            & 0 & 0 \\
0            & 0 & 0
\end{array}
\right), ~
\left( \begin{array}{ccc}
1 & 0 &           0 \\
0 & 1 &           0 \\
0 & 0 &    1   
\end{array}
\right),~
\]
\[
\left( \begin{array}{ccc}
0    & 0            &           0 \\
0    &  1           &           0 \\
0    & \gamma       &           1 
\end{array}
\right), ~
\left( \begin{array}{ccc} 
\varepsilon    & 0            &           0 \\
0              &  1           &           0 \\ 
0              & \gamma       &      1     
\end{array} 
\right),
~{\rm and~}~  
 \left( \begin{array}{ccc}
1            & 0       & 0 \\ 
0            & 1       & 0 \\ 
0            & \gamma  & 1  
\end{array}
\right),  
\]
respectively.  
 Furthermore, it is easy to show that all the 
congruence classes contain at least one compatible matrix, hence, all 
the equivalence classes contain at least one compatible matrix. Thus, the 
number of equivalence classes over ${\bf F}_p$, $p\neq 2$ is $3p+10$.  
\vspace{0.2in} 

\noindent In view of the above discussion, we may now state the following: 

\begin{prop} 
Let $N(s,~1)$ denote the total number of non-isomorphic rings with 
property(T) and characteristic $p$ with maximal Galois subfield ${\bf F}_p$;
and with the same invariants $p$, $n$, $s$, $t$, $\lambda$, where $t=1$.
Then,

\noindent $N(2,~1)= 5 $ or $p+4$ according as $p=2$ or otherwise; 

\noindent $N(3,~1)=11$ or $3p+10$ according as $p=2$ or $p\neq 2$. 

\noindent In general, $N(s,~1) \leq N(s)-1$,    
where $N(s)$ is as in Theorem 5.4. Moreover, this bound is reached when 
$p=2$. 
\end{prop} 

\noindent We next consider the matrices $\beta^{-1}C^{T}AC$, where 
$A\in {\bf M}_s(F)$ is symmetric and $F$ is any finite Galois field
$GF(p^r)$.
Theorem 5.5 gives the number and representatives of congruence classes of 
$s\times s$ symmetric matrices of any positive rank $r\leq s$. Now, if 
$p\neq 2$, for any $s>1$, the number of non-zero congruence classes is 
$2s$. It is easy to verify  that each of these classes contains a 
compatible matrix. Also, as $\beta$ runs over all the elements of 
$F^{*}$, we see that all the  classes of odd rank reduce to one 
equivalence class, namely, to the class with $1$'s in the main diagonal, 
while those of even rank remain distinct.  For instance, 
\[
\left(\begin{array}{ccc} 
1 & 0 & 0 \\
0 & 0 & 0 \\ 
0 & 0 & 0 
\end{array} 
\right) ~ {\rm and~} 
\left(\begin{array}{ccc} 
\epsilon & 0 & 0 \\ 
0        & 0 & 0 \\ 
0        & 0 & 0  
\end{array}  
\right) 
\]
become equivalent to each other. Therefore, the number of
equivalence classes of $1$-dimensional symmetric bilinear forms
over $F$ when $p\neq 2$ is $\frac{3s -1}{2}$ if $s$ is odd
and $\frac{3s}{2}$ if $s$ is even. 
\vspace{0.2in} 

\noindent Now, if $p=2$, the number of non-zero congruence classes of 
$s\times s$ symmetric matrices over $F$ is $\frac{3s -1}{2}$ 
when $s$ is odd; and 
$\frac{3s}{2}$ when $s$ is even. Clearly, as $\beta$ runs over all the 
elements of $F^{*}$, these congruence classes remain  
distinct equivalence classes. Furthermore, it is easy to see that each 
equivalence class contains a compatible matrix; hence, the above numbers 
give the number of non-isomorphic commutative rings with 
characteristic $p$ and of the same invariants. The models of these can 
be deduced from Theorem 5.5. 
\vspace{0.2in} 

\noindent We may then state the following: 
\vspace{0.2in} 

\begin{prop} 
Let $N_c(s,~1)$ denote the total number of isomorphism classes of 
commutative rings with property(T) and characteristic $p$ with maximal 
Galois subfield $GF(p^r)$; and with the same invariants $p$, $n$, $r$, $s$, 
$t$, $\lambda$, where $t=1$. Then 
\[
N_c(s,~1)=\left\{ 
\begin{array}{ll} 
\frac{3s-1}{2}  & ~if~s~is~odd, \\
\frac{3s}{2}    & ~if~s~is~even,   
\end{array} 
\right. 
\]
for any prime $p$. 
\end{prop}      
\vspace{0.2in} 

\noindent In the case where the rings are not commutative and the field $F$ 
is not a prime field, it is intuitively obvious that, in general, there will 
be a lot fewer equivalence classes than congruence classes. Therefore, all 
we can say is that the number of isomorphism classes of the rings in 
question does not exceed $N(s)-1$; the number of non-zero congruence classes. 
\vspace{0.2in} 

\noindent The study of how the congruence classes are subdivided into 
equivalence classes is obviously very important and we  consider 
this in subsequent works.

\subsubsection{The case with $s=2$, $t=2$} 

Let $R$ be a ring with property(T) and characteristic $p$ in which the 
maximal Galois subfield $F$ lies in the centre and with invariants 
$p$, $n$, $r$, $s$, $t$, $\lambda$, where $s=2$ and $t=2$. 
Then, the ring $R$ is 
defined by two structural matrices $A_1$ and $A_2$, where $A_1$ and $A_2$ 
are  $2\times 2$ compatible  matrices over $F$. We know from Lemma 4.1
 how two rings of the same type can be isomorphic with each other.
Moreover, two rings of the same type are isomorphic if and
only if their corresponding spaces of bilinear forms are equivalent.
\vspace{0.2in} 

\noindent Let $N(2,~2)$ denote the number of equivalence classes of 
 2-dimensional spaces of 
$2\times 2$ matrices over $F$ corresponding to  2-dimensional 
spaces of bilinear forms. The number of such equivalence classes may be
determined and the class representatives may be obtained for particular
values of $p$ by using programs we devised that make use of elements from
{\sc MATLAB}. Here we give a representative of such programs in
the Appendix. The number of equivalence classes $N(2,~2)$ is then
given by the following:
  
\[
N(2,~ 2) =\left\{ 
        \begin{array}{ll} 
             10    & {\rm if~} |F|=2 \\ 
             14    & {\rm if~} |F|=3 \\ 
             20    & {\rm if~} |F|=5 \\ 
             26    & {\rm if~} |F|=7.   
        \end{array} 
\right.
\]

\noindent With these results, we may then state the following:  

\begin{prop} 
 The number of mutually non-isomorphic rings with property(T) and  
characteristic $p$ and of the same order with  maximal Galois
subfield ${\bf F}_p$,  and with the same invariants $p$, $n$,
$s$, $t$, $\lambda$, where $s=2$ and $t=2$, is 
\[
\left\{ 
\begin{array}{ll} 
   10     & {\rm if~} p=2, \\
   3p + 5 & {\rm if~} p\neq 2. 
\end{array} 
\right. 
\]
Of these, only  $3$ are commutative (for every prime $p$), the others are not.
\end{prop} 
       
\subsubsection{The case with $s=2$, $t=3$}
 
We now consider the problem of classifying all the 
 rings of a given order with property(T)
and  characteristic $p$, 
in which the maximal Galois subfield $F$ lies in
the centre, for given invariants $p$, $r$, $n$, $s$, $t$, $\lambda$, where
$s=2$ and $t=3$.
\vspace{0.2in} 
 
\noindent Let $R$ be one such ring. Then, the ring $R$ is
defined by three $2\times 2$ compatible structural matrices
$A_1$, $A_2$ and $A_3$ over $F$.
 Then on the basis 
of Lemma 4.1, if $R(D)$ is isomorphic to  
$R(A)$, where $A=\{A_1,~A_2,~A_3\}$ and $D=\{D_1,~D_2,~D_3\}$,
there exist matrices $C$ in
$GL(2,~F)$ and $B=(\beta_{k\rho})$ in $GL(3,~F)$ 
such that 
\begin{eqnarray*} 
D_1 & = & \beta_{11}C^TA_1C +\beta_{12}C^TA_2C +\beta_{13}C^TA_3C, \\
D_2 & = & \beta_{21}C^TA_1C +\beta_{22}C^TA_2C +\beta_{23}C^TA_3C, \\ 
D_3 & = & \beta_{31}C^TA_1C +\beta_{32}C^TA_2C +\beta_{33}C^TA_3C.  
\end{eqnarray*} 
 
\noindent Let $N(2,~3)$ denote the number of equivalence classes of  
3-dimensional spaces of
$2\times 2$ matrices over $F$ corresponding to  3-dimensional
spaces of bilinear forms. Then,  by the {\sc MATLAB} program  
 in the Appendix, we have
\[
N(2,~ 3) =\left\{
        \begin{array}{ll}
             5         & {\rm if~} |F|=2 \\
             7         & {\rm if~} |F|=3 \\  
             9         & {\rm if~} |F|=5  \\ 
\end{array}
\right.
\]
\vspace{0.2in} 
 
\noindent We have partial results to this problem, and
conjecture the number of mutually non-isomorphic rings of order
$p^{nr}$ with property(T) and characteristic $p$.
\vspace{0.2in}

\noindent{\bf Conjecture A}
{\it The number of mutually non-isomorphic rings with property(T) and 
characteristic $p$ and of order $p^{nr}$, 
with  maximal Galois subfield ${\bf F}_p$,  
 and with the same invariants $p$, $n$, $s$, $t$, $\lambda$,
where $s=2$, $t=3$, is
\[
\left\{
\begin{array}{ll}
   5    & {\rm if~} p=2, \\
   p+4  & {\rm if~} p\neq 2.
\end{array}
\right.
\]
Of these, only one is commutative (for every prime $p$), the others are not.}
\vspace{0.1in}

It must be noted that the case $p=2$
follows from  the {\sc MATLAB} program.  
  
\subsubsection{The case with $s=3$, $t=2$} 

By a program similar to that in the Appendix devised using elements from 
{\sc MATLAB}, we find that 
the number of equivalence classes of 2-dimensional spaces 
of $3\times 3$ matrices over 
${\bf F}_2$ corresponding to  2-dimensional spaces of bilinear forms on 3 
variables is 322. All these classes contain at least 
one compatible matrix, and therefore, we conclude that all these classes 
are representatives for the rings in question.
\vspace{0.2in} 

 The number of mutually non-isomorphic rings of characteristic
$2$ with  property(T) may now be given by the following result.

\begin{prop}  The number of mutually non-isomorphic rings
with property(T) and
characteristic $p=2$ with maximal Galois subfield ${\bf F}_2$,  
 and with the same invariants $p$, $n$, $s$, $t$, $\lambda$,
where $s=3$, $t=2$;   is $322$. 
 
 Of these, $14$ are commutative, the others are not.
\end{prop} 

 In the case where the field $F$ is not prime, it is obvious that 
there will be more equivalence classes than we have in the case of 
prime subfields. Therefore, all we can say is that the number of 
isomorphism classes of rings with property(T) and characteristic $p$ 
in which the maximal Galois subring $F$ lies in the centre and with same 
invariants $p$, $n$, $r$, $s$, $t$, $\lambda$, with $s>1$, does not 
exceed  the number of distinct subspaces of ${\bf M}_s(F)$ of 
dimension $t$. 
This upper bound is reached in the case where $t=s^2$, since in this case, 
by Lemma 5.3, we only have one ring for any $s>1$.   
\bigskip

\noindent{\large\bf Acknowledgements}
\vspace{0.2in}

 The author would like to thank Prof. D. Theo for many
comments and illuminating discussions and Prof. B. Choudhary for comments
on an earlier draft of this work.
\bigskip

\noindent{\large\bf References}
\bigskip

\noindent[1] {\bf C. J. Chikunji}, On a Class of Finite Rings;
Communications in Algebra ({\it to appear, Fall 1999}).

\noindent[2] {\bf  P. S. Bremser}, Congruence Classes of Matrices in
$GL_2(F_q)$, Discrete Math. 118 (1993), p243 - 249.

\noindent[3] {\bf W. E. Clark}, A coefficient ring for finite
non-commutative rings, Proc. Amer. Math. Soc. 33, No.1 (1972), p25 - 28.

\noindent[4] {\bf B. Corbas}, Rings with few zero divisors, Math.
Ann. 181 (1969), p1 - 7.

\noindent[5] {\bf B. Corbas \& G. D. Williams}, Matrix Representations
for Three Dimensional Billinear Forms over Finite Fields, Discrete
Mathematics 185 (1998) p51 - 61.

\noindent[6] {\bf M. Hall, JR.}, Combinatorial Theory, Wiley (1986).

\noindent[7] {\bf N. Jacobson}, Structure of rings, Amer. Math. Soc.
Colloq. Publ. XXXVII (1964).

\noindent[8] {\bf M. Newman}, Integral Matrices, Academic Press, N.Y. (1972).

\noindent[9] {\bf R. Raghavendran}, Finite Associative Rings, Compositio
Math. 21, Fasc. 2 (1969), p195 - 229.

\noindent[10] {\bf W. C. Waterhouse}, The number of Congruence Classes in
$M_{n}(F_q)$, Finite Fields and their Applications, 1 (1995), p57 - 63.

\pagebreak

\noindent{\large\bf Appendix}
\vspace{0.2in}

\noindent{\large\bf A {\rm MATLAB} Program for $s=2$, $t=3$ over $F_3$}
{\footnotesize
\begin{verbatim}
function jo(a)
global A
global B
global C
T=[ ];
for i=1:12
  if a >= 2*3^(12 - i)  T(i) = 2;  a = a - 2*3^(12 - i);
   elseif a >= 3^(12 - i)  T(i) = 1; a = a - 3^(12 - i);
   else T(i) = 0;
  end
end
A = [T(1:2); T(3:4)];
B = [T(5:6); T(7:8)];
C = [T(9:10); T(11:12)];

function joh(a)
global M
T = [ ];
for i = 1:4
  if a >= 2*3^(4 - i)   T(i) = 2; a = a - 2*3^(4 - i);
    elseif a >= 3^(4 - i)  T(i) = 1;  a = a - 3^(4 - i);
     else  T(i) = 0;
   end
end
M = [T(1:2); T(3:4)];
 
function john(a)
global N
T = [ ];
for i = 1:9
  if a >= 2*3^(9 - i)   T(i) = 2;  a = a - 2*3^(9 - i);
   elseif  a >= 3^(9 - i)  T(i) = 1;  a = a - 3^(9 - i);
   else  T(i) = 0;
  end
end
N = [T(1:3); T(4:6); T(7:9)];
 
function ph(A, B, C)
global a
a = 3^11*A(1, 1)+3^10*A(1, 2)+3^9*A(2, 1)+3^8*A(2, 2)+
       3^7*B(1, 1)+3^6*B(1, 2)+3^5*B(2, 1)+3^4*B(2, 2)+
         3^3*C(1, 1)+3^2*C(1, 2)+3*C(2, 1)+C(2, 2);
 
x = [1:3^12 - 1];
global x
 
global x
for i = 1:6560   x(i) = 0; end
global A;
global B;
global C;
global M;
global N;
global a;
for i = 6560:3^12 - 1
  jo(i);
    if A == zeros(2)  x(i) = 0;
      if B == zeros(2)  x(i) = 0;
       if C == zeros(2)  x(i) = 0;
         if rem( (A + B), 3) == 0 x(i) = 0;
         if rem( (A + C), 3) == 0 x(i) = 0;
       if rem( (B + C), 3) == 0 x(i) = 0;
      if rem( (A + 2*B), 3) == 0 x(i) = 0;
    if rem( (A + 2*C), 3) == 0 x(i) = 0;
  if rem( (B + 2*C), 3) == 0 x(i) = 0;
if rem( (A + B + C), 3) == 0 x(i) = 0;
  if rem( (A + B + 2*C), 3) == 0 x(i) = 0;
    if rem( (A + 2*B + C), 3) == 0 x(i) = 0;
      if rem( (A + 2*B + 2*C), 3) == 0 x(i) = 0;
      end
    end
  end
end
  end
   end
    end
      end
       end
       end
      end
     end
   end
end
for k = 6561:3^12 - 1
 if x(k) ~= 0
  jo(k);
    for i = 1:80
      joh(i);
        if rem( det(M), 3) ~= 0
          for j = 1:19682
            john(j);
              if rem( det(N), 3) ~= 0
                X = M * A * M';
                Y = M * B * M';
                Z = M * C * M';
                F = N(1, 1)*X + N(1, 2)*Y + N(1, 3)*Z;
                G = N(2, 1)*X + N(2, 2)*Y + N(2, 3)*Z;
                H = N(3, 1)*X + N(3, 2)*Y + N(3, 3)*Z;
                J = rem(F, 3);
                K = rem(G, 3);
                L = rem(H, 3);
                ph(J, K, L);
                if a ~= k  x(i) = 0;
                end
              end
           end
         end
     end
  end
end
 
global x;
n = 0;
  for i = 6561:3^12 - 1
    if  x(i) ~= 0
      n = n + 1;
       jo(i)
       A
       B
       C
     end
  end
n
\end{verbatim}
 }                                                                  
\end{document}